\documentclass[12pt]{article}
\usepackage{epsfig}
\def\C{{\mathbf{C}}}
\def\bC{{\mathbf{\overline{C}}}}

\def\mod{{\mathrm{mod}\,}}
\author{A. Eremenko\thanks{Supported by NSF grant DMS 0100512
and by the Humboldt Foundation}}
\title{Metrics of positive curvature with conic singularities
on the sphere}
\begin{document}
\maketitle
\begin{abstract}
A simple proof is given of the necessary and sufficient condition
on a triple of positive numbers $\theta_1,\theta_2,\theta_3$ for
the existence of
a conformal metric of constant positive curvature on the sphere,
with three conic singularities of total angles $2\pi\theta_1,2\pi\theta_2,
2\pi\theta_3$. The same condition is necessary and sufficient for
the triple $\pi\theta_1,\pi\theta_2,\pi\theta_3$ to be interior angles of
a spherical triangular membrane.
\end{abstract}

The following problem is classical \cite{Picard,Poincare}.
Let $S$ be a compact Riemann surface, $p_1,\ldots,p_n$ points in $S$,
and $\theta_1\ldots,\theta_n$ positive numbers. Does there exist
a conformal Riemannian metric of constant
curvature $K$ with conic singularities
at $p_j$ such that the total angle at $p_j$ is $2\pi\theta_j$?
A complete answer to this question is known when $K\leq 0,$ see
\cite{McO,Troy1,Troy2}.
If $K\leq 0$, a unique metric exists
if and only if this is not prohibited by the Gauss--Bonnet theorem.
In Troyanov's papers \cite{Troy1,Troy2}
the case of non-constant $K$ is also
considered.

On the other hand, little is known for the case $K>0$.
Troyanov \cite{Troy3} considered the case of two points
on the sphere and showed that the necessary and sufficient condition in this
case is $\theta_1=\theta_2$. It follows from the results of Troyanov in
\cite{Troy1} that
for any compact surface $S$, there exists a 
metric of constant positive curvature if the condition
\begin{equation}
\label{00}
0<\chi(S)+\sum(\theta_j-1)<\min\{2,2\min\theta_j\}
\end{equation}
is satisfied.
Here the left inequality comes from the Gauss--Bonnet theorem while the right
one is a technical restriction needed for Troyanov's method to work.
If one assumes that all $\theta_j\in (0,1)$, and $S$ is the sphere,
then condition (\ref{00}) is actually necessary and sufficient \cite{LT}.

In view of this situation, it seems useful to analyse the simple
case of a metric of positive constant curvature on
the sphere with three singularities.
Here we give a complete solution for this case. 
To do this we return to the methods which were used before the
PDE $\Delta u=Ke^u$ was invoked; our arguments are based on the
papers of Riemann \cite{Riemann}\footnote{Apparently the research
of Klein, Schwarz and Poincar\'e which lead to the discovery of
the Uniformization Theorem
was originally motivated by the study of ordinary linear differential
equations, in particular, by the work of Riemann \cite{Riemann}.
Somewhat later, the G\"ottingen Royal Society suggested a
different approach to the uniformization,
a direct study of the non-linear PDE
$\Delta u=e^u$, which was proposed as a topic of a competition.
It seems that the papers \cite{Picard,Poincare} were inspired by this
competition.}. For modern expositions of this work, see,
for example \cite{Ince,Var}.

The results show that for $K>0$ there are complicated restrictions
in addition to the Gauss--Bonnet theorem.

In the generic case, when none of the prescribed angles is a multiple
of $2\pi$,
our result (Theorem 1 below) was earlier obtained by
Umehara and Yamada \cite{Umehara}. The method of these authors is
somewhat indirect (they use a connection with surfaces of constant
mean curvature $1$ in the hyperbolic space.)
and they do not obtain an explicit condition
in the case of integer $\theta_j$.
In a recent preprint \cite{Furuta},
the results which are equivalent to our theorems
1 and 2 are obtained by a pure
geometric method, which is different from the analytic method used here.
It seems that the proofs presented below
are simpler. Conformal metrics of constant curvature on the sphere with
three conic singularities were also studied by physicists  
\cite{Bilal} but they did not address the question of existence or
uniqueness of
such metrics.

We identify $S\backslash\{ p_3\}$ with the complex plane $\C$ and assume
without
loss of generality that $K=1$, $p_1=0,p_2=1$ and $p_3=\infty$.
Suppose that a conformal Riemannian metric of constant curvature $1$ is given on
$D=\C\backslash\{0,1\}$. Then every point $z\in D$ has a neighborhood 
which is isometric to an open set on the unit sphere $\bC$ \cite{Li}.
The isometry
is a conformal map, so we obtain a multi-valued locally univalent
meromorphic function $f$ in $D$, whose monodromy belongs to the group
of orientation-preserving isometries of $\bC$. This $f$ is sometimes
called a
``developing map'' in geometric literature.

The converse is also true: given such a multi-valued function $f$, we define
a conformal
metric of constant curvature $1$ in $D$ by the length element
$$\lambda(z)|dz|=\frac{2|f'(z)||dz|}{1+|f(z)|^2}.$$

The behavior of $f$ at $\{0,1,\infty\}$ reflects the assumption
that these points are conical singularities with total angles
$2\pi\theta_1,2\pi\theta_2,2\pi\theta_3$: in a local coordinate $z$ near
$p_j$ some
fractional-linear transformation $g_j=L_j\circ f$ of $f$ has the form
$$g_j(z)=z^{\theta_j}.$$

It follows that $f$ satisfies a Schwarz differential equation
\begin{equation}
\label{schwarz}
\frac{f'''}{f'}-\frac{3}{2}\left(\frac{f''}{f'}\right)^2=R(z),
\end{equation}
where 
$$R(z)=\frac{1-\theta_1^2}{2z^2}+\frac{\theta_1^2+\theta_2^2-\theta_3^2-1}{2z(1-z)}
+\frac{1-\theta_2^2}{2(z-1)^2}.$$

{\em The equation $(\ref{schwarz})$ defines a metric with desired properties
if and only if its monodromy group is conjugate to a subgroup of
conformal isometries of the sphere}.

It remains to find out when this is the case.
This problem was considered in another context by Arnold and Krylov
in \cite{Ar} but they
state its solution only when $\theta_j\in (0,1)$.

We begin with a remark that 
(\ref{schwarz}) is
an equation with real coefficients. It follows that our function $f$
maps the upper half plane onto a {\em circular triangle} (membrane) $T$
that is a simply connected surface spread
over the sphere whose boundary consists ofthree
arcs of circles (each of these circles 
may be traced more than a full turn). We will show that the monodromy group
is conjugate to a group of conformal isometries if and only if this triangle
$T$ can be mapped by a fractional-linear transformation
onto to a spherical (geodesic) triangle. In the opposite direction,
if a geodesic triangular membrane $T$ is given, we can paste it together
with its mirror image $T'$ and obtain a surface homeomorphic to the sphere,
with metric of positive curvature $1$ and three conic singularities whose total
angles are twice those of $T$.

It is well-known (since the work of Schwarz) that the monodromy group
of the Schwarz equation (\ref{schwarz}) is the same as the projectivized monodromy
group of the hypergeometric (Gauss) equation
\begin{equation}
\label{hyper}
z(1-z)w''+(c-(a+b+1)z)w'-abw=0,
\end{equation}
where 
parameters $a,b,c$ are related to parameters $\theta_j$
in the following way:
\begin{equation}
\label{parameters}
\begin{array}{l}
\pm\theta_1=1-c,\\
\pm\theta_2=a-b,\\
\pm\theta_3=c-a-b.
\end{array}
\end{equation}

1. Let us first consider the case that {\em none of the numbers $\theta_j$
is an integer}.
Let $\Gamma'$ be the group generated by reflections in the sides of
the circular triangle $T$.
Then the subgroup $\Gamma\subset\Gamma'$ of index $2$ consisting of
orientation-preserving transformations coincides with the monodromy group
of the equation (\ref{schwarz}). If $C_1$ and $C_2$ are two sides of $T$ with
the common vertex $v$ then the group $\Gamma$ contains an elliptic
transformation fixing $v$ and the other point $v'$ where the circles containing
$C_1$ and $C_2$ intersect.
It is evident that the group $\Gamma$ is conjugate to a group of isometries
of the sphere if and only if there exists a fractional-linear transformation
which sends all three pairs $(v,v')$ to pairs of diametrically opposite
points. This is the case if and only if $T$ is equivalent by a fractional-linear
transformation to a triangle whose sides are geodesic.

Monodromy of the equation (\ref{hyper})
was explicitly computed by Riemann \cite{Riemann}; this result
is reproduced in many places, for example, in
\cite{Ince,Klein,Picard1,Var}.
It follows from the explicit formulas that (under the condition that none
of the $\theta_j$
is an integer) this monodromy remains
unchanged if a triple $(\theta_1,\theta_2,\theta_3)$ is replaced
by a triple $(\pm\theta_1+m,\;\pm\theta_2+n,\;\pm\theta_3+k)$,
where $(m,n,k)$ are integers with the property $m+n+k\equiv 0\,(\mod 2)$.
Let us call such triples of positive numbers {\em equivalent}.
Every non-integer triple is equivalent to one and only one
triple with the property 
\begin{equation}
\label{3}
0<\theta'_1+\theta'_2\leq 1,\quad 0<\theta'_2+\theta'_3\leq 1,\quad
 0<\theta'_1+
\theta'_3\leq 1.
\end{equation}
It is easy to show that for
every triple of positive numbers satisfying (\ref{3}) there is a unique
(up to a fractional-linear transformation)
circular triangle with angles $(\pi\theta'_1,\pi\theta'_2,\pi\theta'_3)$, and
this triangle consists of one sheet (is a Jordan region on the sphere).
Such triangle $T$ is equivalent to a geodesic triangle if and
only if the sum of its angles is greater than $\pi$.
Thus we obtain
\vspace{.1in}

\noindent
{\bf Theorem 1}. {\em If none of the $\theta_1,\theta_2,\theta_3$ is an integer,
then a conformal metric of constant positive curvature on the sphere with
conic singularities of total angles $2\pi\theta_1,2\pi\theta_2$ and
$2\pi\theta_3$ exists
if and only if the unique equivalent triple with the property
$(\ref{3})$ satisfies $\theta'_1+\theta'_2+\theta'_3>1$.
Such metric of curvature $1$ is unique.}
\vspace{.1in}

The uniqueness statement holds because a triple $(\theta_1,\theta_2,\theta_3)$
uniquely defines the right hand side $R$ of the Schwarz equation (\ref{schwarz}).
A solution of the Schwarz equation with isometric monodromy
is then defined up to a rotation of the sphere.
\vspace{.1in} 

2. Now we consider the case when {\em at least one of the three numbers
$\theta_1,\theta_2,\theta_3$ is an integer}.
We assume that
$$\theta_j\neq 1,\quad j=1,2,3,$$
because the case $\theta_j=1$ is covered by the result of
Troyanov mentioned in the beginning.

First of all, we have a necessary
condition for the monodromy group to be a group of isometries, that
no logarithms are present in the formal solutions of (\ref{hyper})
at its singular points.
(Presence of a logarithm leads to a parabolic
transformation in the monodromy group, see, for example, \cite{Klein}).

Without loss of generality, we consider
the singularity of the equation (\ref{hyper}) at zero
and suppose that $N=\theta_1=1-c$ is
a positive integer. The exponents at zero are $0$ and $N$.
The logarithms will be absent if a power series
\begin{equation}
\label{series}
w=\sum_{k=0}^\infty a_kz^k, \quad a_0\neq 0
\end{equation}
satisfies (\ref{hyper}). Substituting this series to the equation we obtain
$$[k(k-1)+ck]a_k=[(k-1)(k-2)+(a+b+1)(k-1)+ab]a_{k-1}, \quad k=1,2,\ldots.$$
The coefficient in the left hand side is zero when $k=1-c=N$. So a power
series solution of the form (\ref{series})
exists if and only if the coefficient in the right hand side
is zero for some $k=n\in [1,N],$ that is
$$(n-1)(n-2)+(a+b+1)(n-1)+ab=0.$$
or
$$(n-1)(a+b+n-1)+ab=0.$$
Expressing $a,b$ in terms of $\theta_2,\theta_3$ and $N=\theta_1$, from (\ref{parameters}),
we obtain after
simple transformations, that either sum or difference
of $\theta_1$ and $\theta_2$ equals $2n-N-1$. In other words,
one of the equations holds:
\begin{equation}
\label{stra}
|\theta_1\pm\theta_2|\in [0,N-1] \;\mbox{is
an integer of the opposite parity from $N$}.
\end{equation}
This condition is necessary for the monodromy of
(\ref{schwarz}) to consist of isometries.
Let us prove that it is also sufficient.
If the angle at a vertex $v$ of the triangle $T$
is an integer multiple of $\pi$,
and no logarithms are present in solutions, then two sides of $T$ meeting at
$v$ belong to the same circle, so all sides of $T$ belong to two intersecting
circles.
But every such pair of circles can be mapped by a fractional-linear
transformation onto a pair of great circles (by sending their two
points of intersection to diametrically opposite points).
So the group generated by
reflections in sides of $T$ is conjugate to a group of isometries.
Thus
us we obtain
\vspace{.1in}

\noindent
{\bf Theorem 2}. {\em If $\theta_1$ is an integer,
then the necessary and sufficient
condition of the existence of a conformal metric of curvature $1$ on the sphere,
with three conic singularities of angles $2\pi\theta_1,2\pi\theta_2,2\pi\theta_3$
is that
either $\theta_2+\theta_3$ or $|\theta_2-\theta_3|$ is an integer $m$
of opposite parity from $\theta_1$, and $m\leq \theta_1-1$.
This metric
with prescribed conical points is unique when exists.}
\vspace{.1in}

\noindent
{\bf Further remarks}.
\vspace{.1in}

1. If two of the $\theta_j$ are integers, and the condition in Theorem 2
is satisfied, then the third one is an integer as well.
Then the necessary and sufficient condition becomes:

{\em The sum of the three integers is odd, and
each of them is less than the sum of the others.}

The solution of the Schwarz equation in this
case has trivial monodromy and thus it is a rational function.
It has critical points at $0,1$ and $\infty$. The local degrees
at these points are $\theta_j$, and the condition above is necessary and
sufficient for the existence of such rational function.

2. Our theorems 1 and 2 also give necessary and sufficient conditions
for the existence of a spherical triangular membrane with
geodesic sides and prescribed angles. Let us give a precise definition.


A {\em spherical triangular
membrane} is a bordered surface homeomorphic to the closed
unit disc, with Riemannian metric of constant curvature $1$, whose boundary
consists of three geodesic arcs. It is clear how to define an interior angle
of such a membrane. Spherical triangular membranes were considered by
Klein in \cite{Klein} but the following result does not seem
to be explicitly stated there:
\vspace{.1in}

\noindent
{\bf Theroem 3} {\em There exists a spherical triangular membrane with interior
angles $\pi\theta_1,\pi\theta_2,\pi\theta_3$ if and only if the three
numbers $\theta_1,\theta_2$ and $\theta_3$ satisfy the conditions stated in
theorems 1 and 2.}
\vspace{.1in}

It is somewhat hard to visualize spherical triangular membranes which have
angles
greater than $2\pi$; the author thanks Mario Bonk for the following pictures:
\vspace{.1in}

\begin{center}
\begin{picture}(0,0)%
\includegraphics{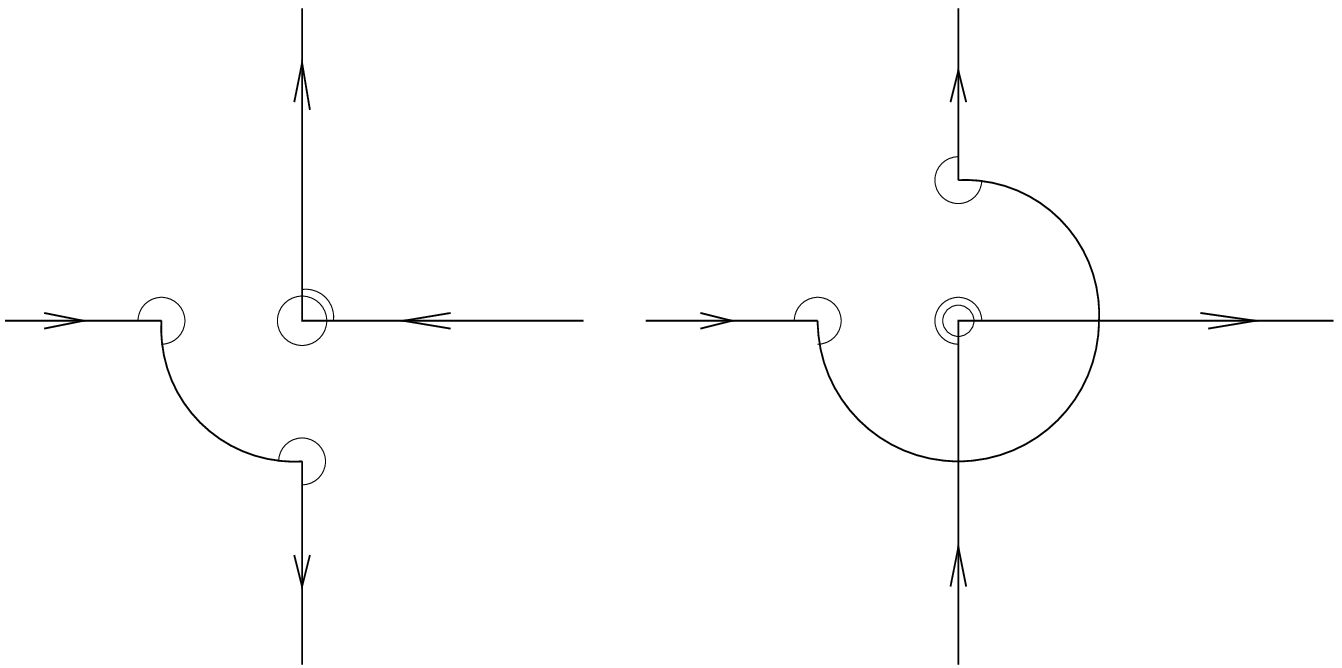}%
\end{picture}%
\setlength{\unitlength}{1973sp}%
\begingroup\makeatletter\ifx\SetFigFont\undefined%
\gdef\SetFigFont#1#2#3#4#5{%
  \reset@font\fontsize{#1}{#2pt}%
  \fontfamily{#3}\fontseries{#4}\fontshape{#5}%
  \selectfont}%
\fi\endgroup%
\begin{picture}(12794,6344)(129,-7283)
\put(1276,-2611){\makebox(0,0)[lb]{\smash{\SetFigFont{12}{14.4}{\rmdefault}{\mddefault}{\updefault}
$1$
}}}
\put(3826,-2686){\makebox(0,0)[lb]{\smash{\SetFigFont{12}{14.4}{\rmdefault}{\mddefault}{\updefault}
$2$
}}}
\put(7651,-2686){\makebox(0,0)[lb]{\smash{\SetFigFont{12}{14.4}{\rmdefault}{\mddefault}{\updefault}
$2$
}}}
\put(9676,-4636){\makebox(0,0)[lb]{\smash{\SetFigFont{12}{14.4}{\rmdefault}{\mddefault}{\updefault}
$1$
}}}
\put(10951,-6286){\makebox(0,0)[lb]{\smash{\SetFigFont{12}{14.4}{\rmdefault}{\mddefault}{\updefault}
$0$
}}}
\put(11026,-2761){\makebox(0,0)[lb]{\smash{\SetFigFont{12}{14.4}{\rmdefault}{\mddefault}{\updefault}
$1$
}}}
\put(1276,-6361){\makebox(0,0)[lb]{\smash{\SetFigFont{12}{14.4}{\rmdefault}{\mddefault}{\updefault}
$0$
}}}
\put(7501,-6286){\makebox(0,0)[lb]{\smash{\SetFigFont{12}{14.4}{\rmdefault}{\mddefault}{\updefault}
$1$
}}}
\end{picture}

\end{center}
\vspace{.1in}

\noindent
{\small {\bf Figure 1.} Spherical triangular membranes with
angles $3\pi/2,3\pi/2,5\pi/2$ (left)
and
\newline
$3\pi/2,3\pi/2,7\pi/2$ (right) represented as surfaces spread over 
the sphere. The sides project into the great circles which are
the coordinate axes and the unit circle.
The covering numbers are shown for each region to help visualize
the surface.} 
\vspace{.2in}

3. In the case of four or more conic singularities solution of the problem
may not be unique.
Indeed, consider the case when all total angles at the singularities
are equal to
$4\pi$. Then solutions $f$ of the Schwarz equation are rational functions,
and the problem is equivalent to finding rational functions with prescribed
simple critical points. Let us call two rational functions
equivalent if they are obtained from each other by post-composition
with a fractional-linear transformation.
It is known \cite{Goldberg}
that for every $2d-2$ points in general position there
exist 
$$u_d=\frac{1}{d}\left(\begin{array}{c}2d-2\\d-1\end{array}\right),\quad
\mbox{the $d$-th Catalan number,}
$$
of equivalence
classes of rational functions of degree $d$ with these critical points.
Thus there is $u_d$ of different conformal metrics of curvature $1$
with conical singularities of total angle $4\pi$ at $2d-2$ given generic points.
The smallest case of non-uniqueness occurs when $d=3$, so there are four
conic singularities with total angles $4\pi$.

Furthermore, for every odd $d\geq 3$,
there are configurations of $2d-2$ points symmetric with respect
to the real line such that the problem 
(with prescribed total angles of $4\pi$ at each of these points)
has no symmetric solutions with
respect to the real line \cite{EG}. On the other hand, if all singularities
belong to the real line, then all solutions of this problem are
symmetric with respect to the real line \cite{EG0}.

{\em Purdue University

West Lafayette, Indiana

eremenko@math.purdue.edu}

\begin{thebibliography}{11}
\bibitem{Ar} V. U. Arnold and Krylov, Uniform distribution of points on
the sphere and some ergodic properties of solutions of ordinary differential
equations in the complex domain, Doklady AN SSSR, N 7 (1963) 47--53.
\bibitem{EG0} A. Eremenko and A. Gabrielov, Rational functions with real
critical points and the B. and M. Shapiro conjecture in real enumerative
geometry, Ann. Math., 155 (2002) 105-129.
\bibitem{EG} A. Eremenko and A. Gabrielov, Counterexamples to pole placement
by static output feedback, to appear in Linear Algebra and Appl., 2002.
\bibitem{Bilal} A. Bilal and J-L. Gervais, Construction of
constant curvature punctured Riemann surfaces with particle-scattering interpretation, J. Geom. Phys., 5 (1988) 277--304.
\bibitem{Goldberg} L. Goldberg, Catalan numbers and branched coverings of the
sphere, Adv. Math., 85 (1991) 129--144.
\bibitem{Furuta} M. Furuta and Y. Hattori, $2$-dimensional singular
spherical space forms, manuscript.
\bibitem{Ince} E. Ince, Ordinary differential equations, Dover, NY, 1956.
\bibitem{Klein} F. Klein, Vorlesungen \"uber die hypergeometrische
Funktionen, Springer, Berlin, 1933.
\bibitem{Li} J. Liouville, Sur l\'equation aux d\'eriv\'ees partielles
$\partial^2\log\lambda/\partial u\partial v\pm2\lambda a^2=0$, J. de Math.,
18 (1853) 71-72.
\bibitem{LT} F. Luo, and G. Tian, Liouville equation and spherical
convex polytopes, Proc. AMS, 116 (1992) 1119--1129.
\bibitem{McO} R. McOwen, Point singularities and conformal metrics
on Riemann surfaces, Proc. AMS, 103 (1988) 222-224.
\bibitem{Picard1} \'E. Picard, Trait\'e d'Analyse, t. III, Gauthier--Villard,
Paris, 1896.
\bibitem{Picard} \'E. Picard, De l'int\'egration de l'\'equation $\Delta u=e^u$
sur une surface de Riemann ferm\'ee, J. reine angew. Math., 130 (1905) 243--258.
\bibitem{Poincare} H. Poincar\'{e}, Fonctions fuchsiennes et l'\'equation
$\Delta u=e^u$, J. de math. pures et appl., (1898) 137--230. 
\bibitem{Riemann} B. Riemann, Beitrage zur Theorie der durch Gauss'sche Reihe
$F(\alpha,\beta,\gamma,x)$ darstellbaren Funktionen, Ges. Math. Werke, 67--83;
Vorlesungen \"uber die hypergeometrische Reihe, Nachtr\"age, III, 69--94.
US edition: Dover, NY, 1953. 
\bibitem{Troy1} M. Troyanov, Prescribing Curvature on compact surfaces with
conical singularities, Trans. AMS, 324 (1991) 793--821.
\bibitem{Troy2} M. Troyanov, Surfaces euclidiennes \`a singularit\'es coniques,
Enseign. Math., 32 (1986) 79--94.
\bibitem{Troy3} M. Troyanov, Metrics of constant curvature on a sphere
with two conical singularities, Lect. Notes Math., 1410, Springer, NY, 296--308
\bibitem{Umehara} M. Umehara and K. Yamada, Metrics of constant curvature
$1$ with three conical singularities on $2$-sphere,
Illinois J. Math., 44 (2000) 72--94.
\bibitem{Var} V. Varadarajan, 
Meromorphic differential equations. Expos. Math. 9 (1991) 97--188. 

\end{thebibliography}
\end{document}